\documentclass[11pt]{article}
\usepackage{graphicx}
\usepackage{amssymb,latexsym,amsmath}
\usepackage{epstopdf}
\newtheorem{lemma}{Lemma}[section]
\newtheorem{theorem}[lemma]{Theorem}
\newtheorem{prop}[lemma]{Proposition}

\newcommand{\g}{\Gamma}
\newcommand{\pf}{\noindent{\em Proof: }}
\newcommand{\epf}{\hfill\hbox{\rule{3pt}{6pt}}\\}
\newcommand{\forme}[1]{}

\begin{document}

\date{\today}
\title{Distance-regular graphs with valency $k$ having smallest eigenvalue at most $-k/2$}
\author{
{\bf Jack Koolen}\footnote{JHK was partially supported by by the National Natural Science Foundation of China
(No. 11471009). }\\School of Mathematical Sciences\\ University of Science and Technology of China, \\Wen-Tsun Wu Key Laboratory of the Chinese Academy of Sciences, \\
230026, Anhui, PR China \\
e-mail: koolen@ustc.edu.cn
\\
\\
{\bf Zhi Qiao}\\School of Mathematical Sciences\\ University of Science and Technology of China, \\
230026, Anhui, PR China \\
email: gesec@mail.ustc.edu.cn}
\maketitle

\section{Introduction}
In this paper, we study the non-bipartite distance-regular graphs with valency $k$ and having a smallest eigenvalue at most
$-k/2$ (For notations and explanation of the graphs, see next section and \cite{bcn} or \cite{DKT}). There are seven infinite families known, namely
\begin{enumerate}
\item The odd polygons with valency 2;
\item The complete tripartite graphs $K_{t,t,t}$ with valency $2t$ at least 2;
\item The folded $(2D+1)$-cubes with valency $2D+1$ and diameter $D \geq 2$;
\item The Odd graphs with valency $k$ at least 3;
\item The Hamming graphs $H(D, 3)$ with valency $2D$ where $D \geq 2$;
\item The dual polar graphs of type $B_D(2)$ with $D \geq 2$;
\item The dual polar graphs of type $^2A_{2D-1}(2)$ with $D \geq 2$.
\end{enumerate}

First we will show a valency bound for distance-regular graphs with a relatively large, in absolute value, smallest eigenvalue. 
\begin{theorem}\label{valbound}
For any real number $ 1> \alpha > 0$ and  any integer $D \geq 2$, the number of coconnected non-bipartite distance-regular graphs with valency $k$ at least two and diameter $D$, having smallest eigenvalue $\theta_{\min}$ not larger than $-\alpha k$, is finite. 
\end{theorem}

{\bf Remarks.} (i) Note that the regular complete $t$-partite graphs $K_{t \times s}$ ($s,t$ positive integers at least $2$) with valency $k = (t-1)s$ have smallest eigenvalue $-s = -k/(t-1)$. 
\\
(ii) Note that there are infinitely many bipartite distance-regular graphs with diameter 3, for example the point-block incidence graphs of a projective plane of order $q$, where $q$ is a prime power.
For diameter 4 this is also true, for example the Hadamard graphs. \\
(iii)  The second largest eigenvalue for a distance-regular graphs behaves quite differently from its smallest eigenvalue. For example $J(n, t)$ $n \geq 2t \geq 4$, has valency $t(n-t)$ and second largest eigenvalue $(n-t-1)(t-1) -1$. So for fixed $t$ there are infinitely many Johnson graphs $J(n,t)$ with second largest eigenvalue larger then $k/2$.\\
\\
Then we classify the non-bipartite distance-regular graphs with diameter at most 4 with valency $k$ having smallest eigenvalue at most $-k/2$ where for diameter 4 we have also the condition $a_1 \neq 0$.

\begin{theorem}\label{class}
Let $\Gamma$ be a non-bipartite distance-regular graph with diameter $D$ at most $4$ where,  if $D=4$, then $a_1 \neq 0$, and valency $k$ at least $2$, having smallest eigenvalue at most $-k/2$.
Then $\Gamma$ is one of the following graphs:
\begin{enumerate}
\item Diameter equals $1$:
 \begin{enumerate}
 \item The triangle with intersection array $\{ 2; 1\}$
 \end{enumerate}
 \item Diameter equals $2$:
 \begin{enumerate}
\item The pentagon with intersection array $\{2, 1; 1, 1\}$;
\item The Petersen graph with intersection array $\{3, 2; 1, 1\}$;
\item The folded $5$-cube with intersection array $\{5, 4; 1, 2\}$;
\item The $3 \times 3$-grid with intersection array $\{ 4, 2; 1, 2\}$;
\item The generalized quadrangle $GQ(2,2)$ with intersection array $\{6, 4; 1, 3\}$;
\item The generalized quadrangle $GQ(2,4)$ with intersection array $\{10, 8; 1, 5\}$;
\item A complete tripartite graph $K_{t,t,t}$ with $t \geq 2$, with intersection array $\{2t, t-1; 1, 
2t\}$;
\end{enumerate}
\item Diameter equals $3$:
\begin{enumerate}
\item The $7$-gon, with intersection array $\{2, 1, 1; 1, 1, 1\}$;
\item The Odd graph with valency $4$, $O_4$, with intersection array $\{4, 3, 3; 1, 1, 2\}$;
\item The Sylvester graph with intersection array $\{5, 4, 2; 1, 1, 4\}$;
\item The second subconstituent of the Hoffman-Singleton graph with intersection array $\{6, 5, 1; 1, 1, 6\}$;
\item The Perkel graph with intersection array $\{6, 5, 2; 1, 1, 3\}$;
\item The folded $7$-cube with intersection array $\{7, 6, 5; 1, 2, 3\}$;
\item A possible distance-regular graph with intersection array $\{7, 6, 6; 1, 1, 2\}$;
\item A possible distance-regular graph with intersection array $\{8, 7, 5; 1, 1, 4\}$;
\item The truncated Witt graph associated with $M_{23}$(see \cite[Thm 11.4.2]{bcn}) with intersection array $\{15, 14, 12; 1, 1, 9\}$;
\item The coset graph of the truncated binary Golay code with intersection array $\{ 21, 20, 16; 1, 2, 12\}$;
\item The line graph of the Petersen graph with intersection array $\{4, 2, 1; 1, 1, 4\}$;
\item The generalized hexagon $GH(2, 1)$ with intersection array $\{ 4, 2, 2; 1, 1, 2\}$;
\item The Hamming graph $H(3,3)$ with intersection array $\{6, 4, 2; 1, 2, 3\}$;
\item One of the two generalized hexagons $GH(2, 2)$ with intersection array $\{6, 4, 4; 1, 1, 3\}$;
\item One of the two distance-regular graphs with intersection array $\{8, 6, 1; 1, 3, 8\}$ (see \cite[p. 386]{bcn});
\item The regular near hexagon $B_3(2)$ with intersection array $\{14, 12, 8; 1, 3, 7\}$;
\item The generalized hexagon $GH(2, 8)$ with intersection array $\{18, 16, 16; 1, 1, 9\}$;
\item The regular near hexagon on $729$ vertices related to the extended ternary Golay code
with intersection array $\{24, 22, 20; 1, 2, 12\}$;
\item The Witt graph associated to $M_{24}$ (see \cite[Thm 11.4.1]{bcn}) with intersection array $\{30, 28, 24; 1, 3, 15\}$;
\item The regular near hexagon $^2A_{5}(2)$ with intersection array $\{42, 40, 32; 1, 5, 21\}$.
\end{enumerate}
\item Diameter equals $4$ and $a_1 \neq 0$;
\begin{enumerate}
\item The generalized octagon $GO(2, 1)$ with intersection array $\{ 4, 2, 2, 2; 1, 1, 1, 2\}$;
\item The distance-regular graph with intersection array $\{6, 4, 2, 1; 1, 1, 4, 6\}$ (see \cite[Thm 13.2.1]{bcn});
\item The Hamming graph $H(4,3)$ with intersection array $\{8, 6, 4, 2; 1, 2, 3, 4\}$;
\item A generalized octagon $GO(2, 4)$ with intersection array $\{10, 8, 8, 8; 1, 1, 1, 5\}$;
\item The Cohen-Tits regular near octagon associated with the Hall-Janko group (see \cite[Thm 13.6.1]{bcn}) with intersection array 
$\{ 10, 8, 8, 2; 1, 1, 4, 5\}$.
\item The regular near hexagon $B_4(2)$ with intersection array $\{30, 28, 24, 16; 1, 3, 7, 15\}$;
\item The regular near hexagon $^2A_{7}(2)$ with intersection array $\{170, 168, 160, 128; 1, 5, 21, 85\}$.
\end{enumerate}
\end{enumerate}
\end{theorem}

 {\bf Remark.} It is not known whether the generalized octagon $GO(2, 4)$ with intersection array $\{10, 8, 8, 8; 1, 1, 1, 5\}$ is unique.

 This result is an extension of De Bruyn's results \cite[Sects. 3.5 \& 3.6]{B06} on regular near hexagons and octagons, with lines with size 3, see also Theorem \ref{rnp}.

As a consequence of Theorem \ref{class}, we also obtain a complete classification of the 3-chromatic distance-regular graphs with diameter 3 and the 3-chromatic distance-regular graphs with diameter 4 and intersection number $a_1 \neq 0$. 

\begin{theorem}\label{3chrom}
(i) Let $\Gamma$ be a $3$-chromatic distance-regular graph with diameter $3$. Then $\Gamma$ is one of the following:
\begin{enumerate}
\item The $7$-gon, with intersection array $\{2, 1, 1; 1, 1, 1\}$;
\item The Odd graph with valency $4$, $O_4$, with intersection array $\{4, 3, 3; 1, 1, 2\}$;
\item The Perkel graph with intersection array $\{6, 5, 2; 1, 1, 3\}$;
\item The generalized hexagon $GH(2, 1)$ with intersection array $\{ 4, 2, 2; 1, 1, 2\}$;
\item The Hamming graph $H(3,3)$ with intersection array $\{6, 4, 2; 1, 2, 3\}$;
\item The regular near hexagon on $729$ vertices related to the extended ternary Golay code
with intersection array $\{24, 22, 20; 1, 2, 12\}$.
\end{enumerate}

(ii) Let $\Gamma$ be a $3$-chromatic distance-regular graph with diameter $4$ and $a_1 \neq 0$. Then $\Gamma$ is the Hamming graph $H(4, 3)$ with intersection array $\{ 8, 6, 4, 2; 1, 2, 3, 4\}$, or the generalized hexagon $GO(2, 1)$ with intersection array $\{ 4, 2, 2, 2; 1, 1, 1, 2\}$.
\end{theorem}

This result is an extension of Blokhuis et al. \cite{BBH07}. In that paper, they determined all the 3-chromatic distance-regular graphs among the known examples. 

This paper is organised as follows, in Section 2 we give definitions and preliminaries, and in Section 3 we give the proof of the valency bound, Theorem \ref{valbound}. In Section 4, we treat the strongly regular graphs. In Section 5 we give a bound on the intersection number $c_2$.
In Sections 6 we treat the case $a_1=1$ and in Section 7 we treat the case $a_1 =0$. In Section 8, we give the proofs of Theorems 1.2 and 1.3.In the last section we give some open problems.

\section{Preliminaries and definitions}
All graphs considered in this paper are finite, undirected and simple (for 
more background information,
see \cite{bcn} or \cite{DKT}). For a connected graph $\Gamma=(V(\Gamma),E(\Gamma))$, the distance $d(x,y)
$ between any two vertices $x,y$ is the length of a shortest
path between $x$ and $y$ in $\Gamma$, and the diameter $D$ is the
maximum distance between any two vertices of $\Gamma$. For any vertex $x$, let $\Gamma_i(x)$ be
the set of vertices in $\Gamma$ at distance precisely $i$ from $x$, where $0 \leq i\leq D$.
For a set of vertices $x_1,\ldots,x_n$, let $\Gamma_1(x_1,\ldots,x_n)$ denote $\cap_{i=1}^n \Gamma_1(x_i)$.
For a non-empty subset $S\subseteq V(\Gamma)$, $\langle S\rangle$ denotes the induced subgraph on $S$.
A graph is coconnected if its complement is connected.
A connected graph $\Gamma$ with diameter $D$ is called a {\em distance-regular graph}~if there are integers
$b_i$, $c_i$ ($0\leq i \leq D$) such that for any
two vertices $x,y\in V(\Gamma)$ with $i=d(x,y)$, there are
exactly $c_i$ neighbors of $y$ in $\Gamma_{i-1}(x)$ and $b_i$
neighbors of $y$ in $\Gamma_{i+1}(x)$. The numbers $b_i, c_i$ and
$a_i:=b_0-b_i-c_i~(0 \leq i\leq D)$ are called the {\em
intersection numbers} of $\Gamma$. Set $c_0=b_{D}=0$. We observe $a_0=0$ and $c_1=1$. The array $\iota(\Gamma)=\{b_0,b_1,\ldots,b_{D-1};c_1,c_2,\ldots,c_{D}\}$
is called the {\em intersection array} of $\Gamma$. In particular $\Gamma$ is a regular graph with valency $k:=b_0$. We define $k_i:=|\Gamma_i(x)|$
for any vertex $x$ and $i=0,1,\ldots,D$. Then we have $k_0=1,~k_1=k,~c_{i+1}k_{i+1}=b_{i} k_{i}~~(0 \leq i\leq D-1)$ and thus
\begin{equation} \label{ki}
k_i=\frac{b_1\cdots b_{i-1}}{c_2 \cdots c_i}k~~~(1\leq i\leq D).
\end{equation}
A regular graph $\Gamma$ on $n$ vertices with valency $k$ is called a 
strongly regular graph with parameters $(n, k \lambda, \mu)$ if there
are two non-negative integers $\lambda $ and $\mu $ such that for any two distinct
vertices $x$ and $y$, $|\Gamma_1(x,y)|=\lambda$ if $d(x,y)=1$ and $\mu$ otherwise. A connected non-complete strongly regular graph is just a distance-regular graph with diameter 2.\\
The adjacency matrix $A=A(\Gamma)$ is the $(|V(\Gamma)|\times |V(\Gamma)|)$-matrix with rows and columns
indexed by $V(\Gamma)$, where the $(x,y)$-entry of $A$ is $1$ if $d(x,y)=1$ and $0$
otherwise. The eigenvalues of $\Gamma$ are the eigenvalues of $A=A(\Gamma)$.
It is well-known that a distance-regular graph $\Gamma$ with diameter $D$ has
exactly $D+1$ distinct eigenvalues $k=\theta_0>\theta_1>\cdots >\theta_D$ which are
the eigenvalues of the following tridiagonal matrix
\begin{equation*}\label{mtx-L}
 L_1:= \left\lgroup
 \begin{tabular}{llllllll}
 $0$ & $k$\\
 $c_1$ & $a_1$ & $b_1$\\
 & $c_2$ & $a_2$ & $b_2$\\
 & & . & . & .\\
 & & & $c_i$ & $a_i$ & $b_i$\\
 & & & & . & . & .\\
 & & & & & $c_{D-1}$ & $a_{D-1}$ & $b_{D-1}$\\
& & & & &\makebox{\hspace{.324cm}} &
 $c_{D}$ & $a_{D}$
 \end{tabular}
 \right\rgroup
\end{equation*}
(cf. \cite[p.128]{bcn}). The standard sequence $\{u_i(\theta)\mid 0 \leq i\leq D\}$
corresponding to an eigenvalue $\theta$ is a sequence satisfying the following recurrence relation
\[
c_iu_{i-1}(\theta) + a_iu_i(\theta) + b_iu_{i+1}(\theta) = \theta u_i(\theta) ~~(1 \leq i \leq D)
\]
where $u_0(\theta) = 1$ and $u_1(\theta) = \frac{\theta}{k}$. Then the multiplicity of eigenvalue $\theta$ is given by
\begin{equation}\label{mi}
m (\theta)=\frac{|V(\g)|}{\sum_{i=0}^{D}k_iu_i^2(\theta)}
\end{equation}
which is known as {\em Biggs' formula} (cf.  \cite[Theorem 4.1.4]{bcn}).

Let $\Gamma$ be a distance-regular graph with valency $k$, $n$ vertices and diameter $D$. For $i=0, 1, \ldots, D$, let $A_i$ be the $\{0, 1\}$-
matrix with where $(A_i)_{xy}= 1$ if and only $d(x, y)=i$ for vertices $x, y$ of $\Gamma$. 
Let ${\cal A}$ be the Bose-Mesner algebra of $\Gamma$, i.e. the matrix algebra over the complex numbers generated by $A= A_1$. Then ${\cal A}$ has as basis $\{ A_0 = I, A_1= A, A_2, \ldots A_D\}$. The algebra ${\cal A}$ also has a basis of idempotents $\{E_0= \frac{1}{n}J, E_1, \ldots, E_D\}$. Define the Krein numbers $q_{ij}^{\ell}$
where $0 \leq i, j, \ell \leq D$ by $E_i \circ E_j = \frac{1}{n} \sum_{\ell = 0}^D q_{ij}^{\ell} E_{\ell}$. 
It is known the Krein numbers are non-negative real numbers, see \cite[Prop. 4.1.5]{bcn}. We will also need the absolute bound. 
Let $\theta_0 = k > \theta_1 . \cdots > \theta_D$ be the distinct eigenvalues of $\Gamma$ with respective multiplicities $m_0 =1, m_1, \ldots, m_D$. Then for $0 \leq i, j, \leq D$ we have 
$$ \sum_{\ell \in \{0, \ldots, D\}, \mbox{ such that } q^{\ell}_{ij}  \neq 0} m_{\ell} \leq \begin{cases} \ \ m_i m_j \ \ \mbox{ if } i \neq j\\
 \ \ m_i (m_i +1)/2 \ \ \mbox{ if } i = j.\end{cases}$$ This is called the absolute bound.

For a graph $\Gamma$, a partition $\Pi= \{ P_1, P_2, \dots P_t\}$  of $V(\Gamma)$ is called equitable
if there are constants $\alpha_{ij}$ $( 1 \leq i, j \leq t)$ such that all vertices $x \in P_i$ have exactly $\alpha_{ij}$ neighbours in $P_j$. The $\alpha_{ij}$'s $(1 \leq i, j \leq t$) are called the parameters of the equitable partition. 

Let $\Gamma$ be a distance-regular graph. For a set $S$ of vertices of $\Gamma$, define 
$S_i := \{ x \in V(\Gamma) \mid d(x, S):= \min\{ d(x,y) \mid y \in S\} =i\}$.
The number $\rho = \rho(S) := \max\{ i \mid S_i \neq \emptyset \}$ is called the covering radius of $S$. 
The set $S$ is called a completely regular code of $\Gamma$ if the distance-partition 
$\{ S= S_0, S_1, \ldots, S_{\rho(S)}\}$ is equitable.
The following result was first shown by Delsarte \cite{D73} for strongly regular graphs and extended by 
Godsil to the class of distance-regular graphs. 

\begin{lemma} \label{delbound} {\em (Delsarte-Godsil bound)}
Let $\Gamma$ be a distance-regular graph with valency $k \geq 2$, diameter $D \geq 2$ and 
smallest  eigenvalue $\theta_{\min}$. 
Let $C \subseteq V(\Gamma)$ be a clique with $c$ vertices.
Then $$c \leq 1 + \frac{k}{-\theta_{\min}},$$
with equality if and only if $C$ is a completely regular code with covering radius $D-1$.
\end{lemma}
A clique $C$ with $\# C = 1 + \frac{k}{-\theta_{\min}}$, is called a Delsarte clique of $\Gamma
$.
It is known that parameters of a Delsarte clique as a completely regular code only depend on 
the parameters of $\Gamma$.

A distance-regular graph $\Gamma$ is called a geometric distance-regular graph if 
$\Gamma$ contains a set of Delsarte cliques $\mathcal C$, such that every edge 
of $\Gamma$ lies in exactly one member $C$ of $\mathcal C$. 

Examples of geometric distance-regular graphs are for example the bipartite distance-regular 
graphs, the Johnson graphs, the Grassmann graphs, the Hamming graphs and the bilinear 
forms graphs. See \cite{KB10}, for more information on geometric distance-regular graphs.

A geometric distance-regular graph with valency $k$ and diameter $D$ is called a regular near $2D$-gon if $c_i a_1 = 
a_i $ for $i=1, 2, \ldots D$. A generalized $2D$-gon of order $(s,t)$, where $s, t \geq 1$ are integers, is a regular 
near $(2D)$-gon with valency 
$k = s(t+1)$ and intersection number $c_{D-1} = 1$. A generalized $4$-gon of order $(s,t)$ is called a generalized quadrangle of order $(s, t)$ and is denoted by GQ$(s,t)$. In similar fashion, a generalized $6$-gon (respectively $8$-gon) of order $(s,t)$ is called a generalized hexagon (octagon) of order $(s, t)$ and is denoted by GH$(s,t)$ (GO$(s,t)$).

\section{Proof of Theorem 1.1}
In this section, we give a proof of the valency bound Theorem \ref{valbound}.
\\
{\bf Proof of Theorem \ref{valbound}:}\\
As $\Gamma$ is coconnected, $\Gamma$ is not complete multipartite. Let $m$ be the multiplicity of $\theta_{\min}$. As $\Gamma$ is not complete multipartite, we see that $k \leq \frac{(m-1)(m+2)}{2}$ holds, by \cite[Thm 5.3.2]{bcn}.

We consider the standard sequence $u_0=1, u_1, \ldots, u_D$ of $\theta = \theta_{\min}$.
Then $u_1 = \theta/k$ and $$u_{i+1} = \frac{(\theta -a_i) u_i -c_i u_{i-1}}{b_i} \ \ (i = 1,2, \ldots, D-1),$$
and as $\theta$ is the smallest eigenvalue $(-1)^iu_i > 0$ for $i=0, 1, \ldots D$, see \cite[Cor. 4.1.2]{bcn}
We may assume that $k \geq 4\alpha^{-2}$, and hence $c_1 =1 \leq (\alpha^{2}/4)k$.
Let $\ell := \max\{ i \mid 1 \leq i \leq D \mbox{ such that } c_i  \leq \alpha^{i+1}2^{-i-1}k \}.$ and let $p := \min\{\ell+1, D\}$.
\\
\\
{\bf Claim 1}
\\
The number $u_i$ satisfies $|u_i| \alpha^{-i} \geq 2^{-i}$ for $i = 0, 1, 2, \ldots, p$.  
\\
{\bf Proof of Claim 1.}\\
We will show it by induction on $i$. 
For $i=0$, it is obvious, as $u_0 =1$. For $i=1$ we have $|u_1| = |\theta/k| \geq \alpha$, so the claim holds for $i=1$.
 
Let $i \geq 2$. Then $u_i = \frac{(\theta-a_{i-1})u_{i-1} - c_{i-1} u_{i-2}}{b_{i-1}}$ holds.
As $(-1)^i u_i > 0, c_{i-1} \leq \alpha^i 2^{-i}k, b_{i-1} \leq k, a_{i-1} \geq 0, |u_{i-2}| \leq 1, \theta \leq -\alpha k$ we see that $$|u_i| \geq \frac{ k  \alpha |u_{i-1}| - \alpha^i2^{-i}k}{k}.$$
By the induction hypothesis we obtain $|u_i| \alpha^{-i} \geq 2^{1-i} -  2^{-i}= 2^{-i}$.
This shows the claim by induction.
\epf

Let $1 \leq q \leq p$ be such that $k_q$ is maximal among $k_0, k_1, \ldots, k_p$. 
\\
\\
{\bf Claim 2}\\
The number of vertices $n$ of $\Gamma$ satisfies $n \leq (D+1) 2^{(q+1)(D-q)}\alpha^{(q+1)(q-D)}k_q$.
\\
{\bf Proof of Claim 2.}\\
If $q < p$, or if $q=D$, then $n = \sum_{i=0}^D k_i \leq (D+1)k_q$, as, then $k_q = \max\{ k_i \mid 0 \leq i \leq D\}$ and hence the claim follows in this case.\\
So we may assume $p=q <D$. 
As $c_q > \alpha^{q+1}2^{-(q+1)}k$ and 
$$ k_{q+j} = k_q \frac{b_q b_{q+1} \ldots b_{q+j-1}}{c_{q+1} c_{q+2} \ldots c_{q+j}} < k_q k^j c_q^{-j} < 
k_q 2^{(q+1)j}\alpha^{-(q+1)j}$$ for $j=0, 1, \ldots, D-q$. As $n = \sum_{i=0}^D k_i$, Claim 2 follows.
\epf

Let $f(D, \alpha) := \max\{ (D+1)2^{(q+1)(D-q)+2q}\alpha^{(q+1)(q-D)-2q} \mid q=1,2, \ldots D\}$.

By Biggs' formula, see \cite[Thm 4.1.4]{bcn}, Claim 1 and Claim 2, we have
$$m = \frac{n}{\sum_{i=0}^D u_i^2 k_i }< \frac{n}{u^2_q k_q} \leq (D+1) 2^{(q+1)(D-q)}\alpha^{(q+1)(q-D)} \alpha^{-2q} 2^{2q} \leq f(D, \alpha).$$
Hence $k \leq (f(D, \alpha)-1)(f(D, \alpha) +2)/2.$ This shows the theorem with $\kappa(D, \alpha) = (f(D, \alpha)-1)(f(D, \alpha) +2)/2.$
\epf

\section{Diameter 2}
In this section we will determine the connected strongly regular graphs with valency $k \geq2$ and smallest eigenvalue at most $-k/2$.
\begin{prop}\label{diam=2}
Let $\Gamma$ be a non-complete non-bipartite connected strongly regular graph, valency $k \geq 2$ and 
smallest eigenvalue $\theta_{\min}$ satisfying $\theta_{\min} \leq -k/2$, then $\Gamma$ is 
one of the following:
 \begin{enumerate}
\item The pentagon with intersection array $\{2, 1; 1, 1\}$;
\item The Petersen graph with intersection array $\{3, 2; 1, 1\}$;
\item The folded $5$-cube with intersection array $\{5, 4; 1, 2\}$;
\item The $3 \times 3$-grid with intersection array $\{ 4, 2; 1, 2\}$;
\item The generalized quadrangle $GQ(2,2)$ with intersection array $\{6, 4; 1, 3\}$;
\item The generalized quadrangle $GQ(2,4)$ with intersection array $\{10, 8; 1, 5\}$;
\item A complete tripartite graph $K_{t,t,t}$ with $t \geq 2$, with intersection array $\{2t, t-1; 1, 
2t\}$
\end{enumerate}
\end{prop}

Before we show this proposition we recall the following classification of Seidel. \begin{theorem}{\em (Seidel \cite{S68}, see also \cite[Thm 3.12.4(i)]{bcn})}\label{seidel}
Let $\Gamma$ be a strongly regular graph and with second smallest eigenvalue $-2$. 
Then $\Gamma$ is one of the following graphs:
\begin{enumerate}
\item A Cocktail Party graph $K_{n \times 2}$, with $n \geq 2$;
\item A $t \times t$-grid with $t \geq 2$;
\item A triangular graph $T(n)$ with $n \geq 4$;
\item The Petersen graph;
\item The Schl\"{a}fli graph;
\item The Shrikhande graph;
\item One of the three Chang graphs;
\item The halved 5-cube
\end{enumerate}
\end{theorem}
From now on let $\Gamma$ be a non-bipartite distance-regular graph with valency $k \geq 
2$, 
diameter $D \geq 2$ and smallest eigenvalue $\theta_{\min} \leq -k/2$. 
If $ a_1 \geq 1$, then, by Lemma \ref{delbound}, $\Gamma$ has no 4-cliques and any triangle is a completely 
regular code.  \\
\\
Let us first consider the case $a_1 \geq 2$. Then any triangle $T= \{x, y, z\}$ is a completely
regular code and any vertex $u$ at distance 1 from $T$ has exactly two neighbours in $T$.
Let $A_{ab} := \{ u \in V(\Gamma) \mid u \sim a, u \sim b\}$ where $a \neq b$ and $a, b \in 
\{x, y, z \}$. Then $A_{ab}$ forms a coclique, as there are no 4-cliques by the Delsarte-Godsil bound, and if $\{a, b, c\} = \{x, y, z\}$, then each vertex of $A_{ab}$ is adjacent to each vertex of $A_{ac}$.
As the valency of $x$, $y$ and $z$  equals $\# A_{xy} + \# A_{xz}$,  $\# A_{xy} + \# A_{yz}$,  
and $\# A_{xz} + \# A_{yz}$, respectively, it follows that $\Gamma$ is the complete tripartite 
graph $K_{t, t, t}$ where $t = \# A_{xy} = \# A_{xz} =\# A_{yz}$.
This shows:
\begin{lemma}\label{completetripartite}
Let $\Gamma$  be a distance-regular graph with valency $k \geq 2$, diameter $D \geq 2$ 
and smallest eigenvalue $\theta_{\min}$. If $\theta_{\min} \leq -k/2$ and $a_1 \geq 2$,
then $\Gamma$ is a complete tripartite graph $K_{t,t,t}$ for some $t \geq 2$.
\end{lemma}
This shows that, if the distance-regular graph is coconnected, then $a_1 \leq 1$.

Now we are ready to give the proof of Proposition \ref{diam=2}\\
\\
\pf Assume the graph is not bipartite. 
First let us discuss the case when $\theta_{\min}$ is not an integer.
Then $\Gamma$ has intersection array $\{ 2t, t; 1, t\}$ and smallest eigenvalue 
$\frac{-1-\sqrt{4t+1}}{2}$. Hence $\theta_{\min} \leq -k/2 = -t$ implies that $t \leq 2$, and we have that $\Gamma$ is the pentagon as for $t=2$, $\theta_{\min}$ is an integer.
So from now we may assume that $\theta_{\min}$ is an integer.
 Let $\theta_1$ be the other non-trivial eigenvalue of $\Gamma$. Then  $\theta_1$ is an 
 non-negative integer. It follows that $c_2 -k = \theta_1 \theta_{\min}\leq -k\theta_1/2$. This 
 implies that $\theta_1 \leq1$. For $\theta_1=0$, we obtain the complete tripartite graphs, and 
 for $\theta_1 =1$, the complement of $\Gamma$ has smallest eigenvalue $-2$. These have been classified in Theorem \ref{seidel} and by checking them we obtain the proposition.
 \epf
 
 \section{A bound on $c_2$}
In this section we will give a bound on $c_2$. We first give the following result. This is a slight generalisation of 
\cite[Prop. 4.4.6 (ii)]{bcn}. We give a proof for the convenience for the reader, following the proof of \cite{bcn}.
Before we do this we need to introduce the
following.
 Let $\Gamma$ be a distance-regular graph with valency $k$ with an eigenvalue 
$\theta$, say with multiplicity $m$. Let $1=u_0, u_1, \ldots, u_D$ be the standard sequence of $\theta$.
Then there exists a map $\phi: V(\Gamma) \rightarrow {\bf R}^m: x \mapsto \overline{x}$, such that
the standard inner product between $u$ and $v$ satisfies $\langle \overline{u}, \overline{v} \rangle =  
u_{d_{\Gamma}(u,v)}$.
\begin{prop} \label{propcompletebip}{\em (Cf. \cite[Prop. 4.4.6(ii)]{bcn}}
Let $\Gamma$ be a distance-regular graph with diameter $D \geq 2$, and valency $k \geq 2$.
Assume $\Gamma$ contains an induced $K_{r,s}$ for some positive integers $r$ and $s$. Let $\theta$ 
be an eigenvalue of $\Gamma$, distinct from $\pm k$, with standard sequence $1=u_0, u_1, \ldots, u_D$.
Then 
\begin{equation} \label{eq1} 
(u_1 + u_2) ((r+s)\frac{1-u_2}{u_1 + u_2} +2rs) \geq 0
\end{equation} and 
\begin{equation}\label{eq2} (u_1 - u_2) ((r+s)\frac{1-u_2}{u_1 - u_2} -2rs) \geq 0
\end{equation}hold.
In particular, if $\theta$ is the second largest eigenvalue, then $$1+ \frac{b_1}{\theta +1}\geq \frac{2rs}{r+s}$$ holds. 
\end{prop}
\pf Let $X= \{x_1, x_2, \ldots, x_r\}$ and $Y= \{y_1, y_2, \ldots, y_s\}$ be the two color classes of the 
induced $K_{r,s}$. Let $G$ be the Gram matrix with respect to the set  $\{ \overline{u} \mid u \in X \cup Y\}$.
Now $\{X, Y\}$ is an equitable partition of $G$ with quotient matrix $Q$ where 
$$Q = \begin{pmatrix}1+ (r-1)u_2   & su_1 \\
ru_1  & 1  + (s-1)u_2
\end{pmatrix}.$$
Multiplying the first column of $Q$ by $s$ and the second column by $r$ we obtain 
the matrix $$Q' := \begin{pmatrix}s + s(r-1)u_2   & rsu_1 \\
rsu_1  & r  + r(s-1)u_2
\end{pmatrix}.$$
  
  As $G$ is positive semi-definite, it follows that $Q$ and $Q'$ are both positive semi-definite and hence
  $(1 \ 1) Q' (1 \ 1)^T \geq 0$ and $(1, -1) Q'(1, -1)^T\geq 0$. 
  Hence we obtain Equations (\ref{eq1}) and (\ref{eq2}). 
  If $\theta$ is the second largest eigenvalue of $\Gamma$, then $1 > u_1 > u_2> \cdots > u_D$ (using that the largest eigenvalue  of the matrix $T$ of \cite[p. 130]{bcn} equals $\theta_1$), and $1 + \frac{b_1}{\theta_1} = \frac{1-u_2}{u_1 - u_2}$ both hold. This implies the in particular statement. 
\epf

This leads us to the following result.
\begin{lemma}\label{c2bound}
Let $\Gamma$ be a non-bipartite distance-regular graph with diameter $D \geq 3$ and valency $k \geq 2$. If the smallest eigenvalue of $\Gamma$, $\theta_{\min}$, is at most $-k/2$, then $a_1 \leq 1$ and $c_2 \leq 5 + a_1$. 
\end{lemma}
\pf
Let $\theta:= \theta_{\min} \leq -k/2$. We already have established that if $a_1 \geq 2$, then the graph is 
complete tripartite and hence diameter is equal to 2. So this implies $a_1 \leq 1$. 
Let $1=u_0, u_1, \ldots, u_D$ be the standard sequence of $\Gamma$ with respect to 
$\theta$. Then $u_1 +u_2 = \frac{1}{kb_1}(\theta+k)(\theta+1) <0$.
The induced subgraph of $\Gamma$ consisting of two vertices at distance 2 and their common 
neighbours is a $K_{2,c_2}$. By  Equation (\ref{eq1}), we obtain that 
if $a_1 =0$, then $3 > 1 - \frac{k-1}{\theta_{\min}-1} \geq \frac{4c_2}{2+ c_2}$ and hence $c_2 \leq 5$.
If $a_1 =1$ then $\theta_{\min} = -k/2$ and $\sigma_1 = -1/2$ and $\sigma_2 = 1/4$ and again using  Equation (\ref{eq1}), we obtain $c_2 \leq 6$. This shows the lemma.
\epf

 \section{The case $a_1 =1$}
 
 In this section we will discuss the situation for $a_1=1$.
 We will start with the following easy observation.
 \begin{prop}\label{a1}
 Let $\Gamma$ be a distance-regular graph with valency $k \geq 3$, diameter 
 $D \geq 2$, intersection number  $a_1=1$, and smallest eigenvalue 
 $\theta_{\min} \leq -k/2$. 
 Then $\theta_{\min} = -k/2$, $\Gamma$ is geometric, and there exists an integer $i, 2 \leq i \leq D$,
 such that $a_j = c_j$ for $1 \leq j < i$,  $a_i = k/2$ and $a_j = b_j$ for 
 $i+1 \leq j \leq D$, with the understanding that $b_D = 0$. 
 Moreover, if $a_D = k/2$, then $\Gamma$ is a regular near $2D$-gon of order $(2, k/2-1)$. 
 \end{prop}
 \pf
 As each triangle is a Delsarte clique, it follows that $\Gamma$ is a geometric distance-regular graph. This implies the proposition. (For details, we refer to Koolen and Bang \cite{KB10}.)
 \epf

In the following result, we summarise the known existence results about regular near $2D$-gons with $a_1 =1$.

\begin{theorem}\label{rnp} Let $D \geq 2$.
Let $\Gamma$ be a regular near $2D$-gon of order $(2,t)$.
Then $c_2 \in \{1,2,3,5\}$. Moreover,
the following holds:
\begin{enumerate}
\item If $c_2 = 5$, then $\Gamma$ is the dual polar graph of type $^2A_{2D-1}(2)$.
\item If $c_2 =3$, then $\Gamma$  is the  dual polar graph of type $B_D(2)$, or if $D=3$, the Witt graph associated to $M_{24}$ (see \cite[Thm 11.4.1]{bcn}) with intersection array $\{30, 28, 24; 1, 3, 15\}$;.
\item If $c_2=2$, then $\Gamma$ is the Hamming graph $H(D, 3)$, or if $D=3$, the coset graph of the truncated binary Golay code with intersection array $\{ 21, 20, 16; 1, 2, 12\}$;
\item If $c_2=1$ and $D=3$, then $\Gamma$ is one of the following:
\begin{enumerate}
\item The generalized hexagon $GH(2, 1)$ with intersection array $\{ 4, 2, 2; 1, 1, 2\}$;
\item The two generalized hexagons $GH(2, 2)$ with intersection array $\{6, 4, 4; 1, 1, 3\}$;
\item The generalized hexagon $GH(2, 8)$ with intersection array $\{18, 16, 16; 1, 1, 9\}$.
\end{enumerate}
\item If $c_2=1$ and $D=4$, then $\Gamma$ is one of the following:
\begin{enumerate}
\item The generalized octagon $GO(2, 1)$ with intersection array $\{ 4, 2, 2, 2; 1, 1, 1, 2\}$;
\item A generalized octagon $GO(2, 4)$ with intersection array $\{ 10, 8, 8, 8; 1, 1, 1, 5\}$;
\item The Cohen-Tits regular near octagon associated with the Hall-Janko group (see \cite[Thm 13.6.1]{bcn}) with intersection array 
$\{ 10, 8, 8, 2; 1, 1, 4, 5\}$.
\end{enumerate}
\end{enumerate}
\end{theorem}
\pf 
For $D=2$, see for example \cite[Cor. 10.9.5]{GR01}.
For $D=3$, see \cite[Sect. 3.5]{B06}.
For $D=4$, this follows from \cite[Sect. 3.6]{B06} and \cite{B10}.
This shows the theorem for diameter at most 4. 
The cases $c_2 =5$ and $c_2 =4$ follows from 
Brouwer and Wilbrink  \cite{BW83} who classified the regular near $2D$-gons having
$D \geq 4$, $c_2 \geq 3$ and  $a_1 \geq 1$, see also \cite[Thm 9.11]{DKT}.
For $D \geq 3$, $c_2 \geq 2$, $c_3 =3$ and  $a_1 \geq 1$,
it was shown by Van Dam et al. \cite[Thm 9.11]{DKT} that $\Gamma$ is the Hamming graph $H(D, 3)$. 
For $D \geq 4$ vand $c_2=2$, it is shown by Brouwer and Wilbrink that then also the regular near octagon with intersection array $\{ 2c_4,  2c_4 - 2, 2c_4 - 2c_2, 2c_4 - 2c_3; 1, c_2, c_3, c_4\}$ must exist and as by the diameter 4 case, we find that it must be the Hamming graph $H( 4, 3)$, and hence $c_3 = 3$ holds. This shows that if $c_2=2$ and $D \geq 4$ we must have the Hamming graph $H(D, 3)$. This finishes the proof of the theorem.
\epf

For $a_1=1$ and $\theta_{\min} = -k/2$, we can improve the valency bound of Theorem \ref{valbound}.
Note that in Hiraki and Koolen \cite{HK04} a similar bound was obtained for regular near polygons.
\begin{prop}\label{valbounda1}
Let $\Gamma$ be a distance-regular graph with $a_1$, valency $k\geq 4$ and diameter $D \geq 2$.
Then $k \leq 2^{2D+1}-2$.
Moreover, if $c_D = k$, then $k \leq 2^{2D-2} -2$.
\end{prop}
\pf
Let $(u_0, u_1, \ldots, u_D)$ be the standard sequence corresponding to $\theta_{\min} =-k/2$. Let $m$ be the multiplicity of $\theta_{\min}$. By Proposition \ref{a1} there exists  $2\leq i \leq D$ be such that $a_i=k/2$.  
It is easy to show by induction, again using Proposition \ref{a1}, that $u_j = (-2)^{-j}$ if $j \leq i$ and $u_j = (-2)^{j-2i}$ if $j \geq i$. 

So, by Biggs' formula, $m \leq \frac{1}{\min\{u_i^2\mid i=1,2,\ldots, D\}} \leq 2^{2D}$. As $2m \geq k +2$, we find $k \leq 2^{2D+1} -2.$

If $i \leq D-1$, then $$k_{i+1} +k_{i-1} \geq \frac{b_i + c_i}{\max\{c_{i+1},b_{i-1}\}} k_i \geq \frac{k/2}{k} k_i= k_i/2.$$
For positive real numbers $a,b,c,d$, if $a/c \geq b/d$ holds, then 
$$\frac{a}{c}  \geq \frac{a+b}{c+d} \geq \frac{b}{d}$$ holds. Using this, we see 
that
$$m = \frac{\sum_{j=0}^D k_j}{\sum_{j=0}^D k_j u_j^2} \leq \frac{ k_{i-1} + k_i + k_{i+1}}
{2^{-2i+2} k_{i-1} + 2^{-2i} k_i + k_{i+1}2^{-2i+2}},$$ as 
$|u_j| \geq 2^{2-i}$ for $j \not \in \{ i-1, i, i+1\}$.  
Hence $m \leq \frac{3/2}{3 \times 2^{-2i}}= 2^{2i-1} \leq 2^{2D-3}.$
As $k \leq 2m-2$, by \cite[Prop. 3]{HK04}, we see that $k \leq 2^{2D-2} -2$ if $i \leq D-1$. This shows the proposition as $c_D=k$ if and only if $i \leq D-1$.
\epf

In view of Proposition \ref{a1} and Theorem \ref{rnp}, we only need to classify the distance-regular graphs with diameter $D$ equals to 3 or 4, $a_1 =1$ and $c_D=k$. 

\begin{theorem}\label{a1class}
Let $\Gamma$ be a distance-regular graph with diameter $D$ equals $3$ or $4$, $a_1 =1$ and valency $k$, $c_D = k$ and smallest eigenvalue $-k/2$.
Then one of the following hold:
\begin{enumerate}
\item $D=3$ and $\Gamma$ is a distance-regular graph with intersection array $\{8, 6, 1; 1, 3, 8\}$ (see \cite[p. 224]{bcn}), or the line graph of the Petersen graph with intersection array $\{ 4, 2, 1; 1, 1, 4\}$;
\item $D=4$ and $\Gamma$ is the distance-regular graph with intersection array $\{6, 4, 2, 1; 1, 1, 4, 6\}$ (see \cite[Thm 13.2.1]{bcn}).
\end{enumerate}
\end{theorem}
\pf As $a_1 =1$, the valency $k$ is even. By Proposition \ref{valbounda1} we have for diameter 3 that the valency $k$ is bounded by $k \leq 14$ and 
for diameter 4 we obtain $k \leq 62$. We generated all the possible intersection arrays of diameter 3 and 4 with 
$k \leq 14$ for diameter 3 and $k \leq 62$  for diameter 4, such that $c_2 \leq 6$, the $c_i$'s are increasing, the $b_i$'s are
decreasing, the valencies $k_i$ are positive integers, satisfying the conditions of Proposition \ref{a1}, $c_d =k$ and  the multiplicities of the eigenvalues are positive integers. 
Besides the intersection arrays in the theorem we obtained  only the following two intersection arrays $\{ 10, 8, 3; 1, 2, 10\}$ and
$\{12, 10, 3; 1, 3, 12\}$. As both have eigenvalue $-k/2$ with multiplicity $7$ and the number of vertices equals 63, we find by the absolute bound (see \cite[Prop. 4.1.5]{bcn}) that if the graph exists, it must have a vanishing Krein parameter, but that is not the case. So there is no distance-regular graph with  either of these two intersection arrays. 
This shows the theorem.
\epf

 \section{Diameter 3 and $a_1 =0$}
Note that there are infinitely many bipartite distance-regular graphs with diameter 3. 

In the following result, we show that a non-bipartite distance-regular graph with 
diameter 3, valency $k$ and smallest eigenvalue at most $-k/2$ has $k$ at 
most 64.

\begin{prop}\label{bounda1=0}
Let $\Gamma$ be a non-bipartite triangle-free distance-regular graph with diameter $3$, valency $k$ and smallest eigenvalue $\theta_{\min}$ at most $-k/2$.
Then $k \leq 64$ holds.
\end{prop}
\pf
Let $(u_0=1, u_1, u_2, u_3)$ be the standard sequence with respect to $
\theta:= \theta_{\min}$.  
Let $\Gamma$ have distinct eigenvalues $\theta_0= k > \theta_1 > \theta_2 > \theta_3 =\theta$. 
Let $$L = \begin{pmatrix} 0  & k & 0 & 0 \\
1 & 0 & k-1 & 0\\
0 & c_2 & a_2 & b_2\\
0 & 0 & c_3 & a_3
\end{pmatrix}.$$
The matrix $L$ has as eigenvalues $\theta_0, \theta_1, \theta_2, \theta_3$, and hence tr$(L^2) = \theta_0^2 + \theta_1^2 + \theta_2^2 + \theta_3^2 \geq 
k^2 +k^2/4$. 
On the other hand we have tr$(L^2) = a_2^2 + a_3^2 + 2k + 2c_2(k-1) + 2c_3 b_2$. Replacing $b_2$ by $k-a_2 -c_2$ and $a_3$ by $k- c_3$, we obtain
tr$(L^2) = k^2 + (c_3 -a_2)^2 + k(2+2c_2) -2c_2(1+c_3) \leq k^2 + (a_2 -c_3)^2
+12k$ as $c_2 \leq 5$ by Lemma \ref{c2bound}. This means that $(a_2 -c_3)^2 + 12k \geq k^2/4$. Now assume $k \geq 65$ to obtain a contradiction. Then $12k \leq (12/65)k^2$. 
This implies that $|a_2 - c_3| \geq (0.255)k$. This means that at least one of $c_3$ and $a_2$ is at least $(0.255)k$. 
We are going to estimate the multiplicity $m$ of $\theta$. 
On the one hand, $m \geq k$, as $a_1 =0$. On the other hand,
by Biggs' formula, we have $$m = \frac{n}{\sum_{i=0}^3 k_i u_i^2},$$
where $n$ is the number of vertices of $\Gamma$. 
This means that $$m \geq\frac{1}{\min\{ u_i^2 \mid i=0,1,2,3\}}.$$
We have 
$u_0 =1$, $u_1 = \theta/k \leq -1/2$, $$u_2 = \frac{\theta u_1 -1}{k-1} \geq \frac{k-4}{4k-4} \geq \frac{61}{256}.$$ So in order that $m \geq k \geq 65$ holds, we must have  $u_3^2 < 1/64$, or, in other words, $u_3 > -1/8$. 
We obtain $$(k-a_2)/8> b_2/8>-b_2 u_3 = (-\theta+a_2)u_2 +c_2 u_1 \geq (k/2 +a_2)(61/256)
-5, $$
and this implies $a_2 < \frac{3k + 2560}{186}$. As $k \geq 65$, it follows that 
$a_2 < k/4$. As we have already established that at least one of $a_2$ and $c_3$ is at least $(0.255)k$, we find $c_3 \geq (0.255)k$. 
We find $-b_2 u_3 = u_2 ( a_2 - \theta) + c_2 u_1 \geq \frac{61}{256}(-\theta) + 5\frac{\theta}{k}= -\theta(\frac{61}{256} -\frac{5}{k})\geq \frac{k}{2}(\frac{61}{256} -\frac{5}{k})$. As $b_2 < k$ and $k \geq 65$, we obtain
$-u_3 \geq \frac{1}{2}(\frac{61}{256} -\frac{5}{65}) > 41/512$.
  
For positive real numbers $a,b,c,d$, if $a/c \geq b/d$ holds, then 
$$\frac{a}{c}  \geq \frac{a+b}{c+d} \geq \frac{b}{d}$$ holds. Using this, we see 
that
$$m = \frac{\sum_{i=0}^3 k_i}{\sum_{i=0}^3 k_i u_i^2} \leq \frac{k_2 + k_3}
{u^2_2 k_2 + u^2_3k_3}$$ holds, as $\frac{1+k}{1 + k u_1^2} \leq 4$.
Now $$\frac{k_2 +k_3}{k_2 u_2^2 + k_3 u_3^2} = 
\frac{ c_3 +b_2}{c_3 u_2^2+b_2 u_3^2}\leq \frac{ c_3 + k}{c_3 u_2^2 + k u_3^2},$$ as $|u_3| < 1/8 < u_2$, $k_3 = \frac{b_2}{c_3}k_2$ and $b_2 < k$ all hold.
Using $c_3 \geq (0.255)k$, $u_2 \geq 61/256$ and $|u_3| \geq 41/512$, we find that $$m \leq \frac{k_2 +k_3}{k_2 u_2^2 + k_3 u_3^2} \leq 64,$$
a contradiction. This shows the proposition.\epf

Now we come to the main result of this section.
\begin{theorem}\label{classa=0}
Let $\Gamma$ be a non-bipartite distance-regular graph with diameter 3, $a_1=0$, valency $k \geq 2$ and smallest eigenvalue at most $-k/2$. Then $\Gamma$ is one of the following:
\begin{enumerate}
\item The $7$-gon, with intersection array $\{2, 1, 1; 1, 1, 1\}$;
\item The Odd graph with valency $4$, $O_4$, with intersection array $\{4, 3, 3; 1, 1, 2\}$;
\item The Sylvester graph with intersection array $\{5, 4, 2; 1, 1, 4\}$;
\item The second subconstituent of the Hoffman-Singleton graph with intersection array $\{6, 5, 1; 1, 1, 6\}$;
\item The Perkel graph with intersection array $\{6, 5, 2; 1, 1, 3\}$;
\item The folded $7$-cube with intersection array $\{7, 6, 5; 1, 2, 3\}$;
\item A possible distance-regular graph with intersection array $\{7, 6, 6; 1, 1, 2\}$;
\item A possible distance-regular graph with intersection array $\{8, 7, 5; 1, 1, 4\}$;
\item The truncated Witt graph associated with $M_{23}$(see \cite[Thm 11.4.2]{bcn}) with intersection array $\{15, 14, 12; 1, 1, 9\}$;
\item The coset graph of the truncated binary Golay code with intersection array $\{ 21, 20, 16; 1, 2, 12\}$;
\end{enumerate}
\end{theorem}
\pf By Proposition \ref{bounda1=0} we have that the valency $k$ is bounded by $k \leq 64$. We generated all the possible intersection arrays of diameter 3 with 
$k \leq 64$, such that $a_1 =0$, $c_2 \leq 5$, the $c_i$'s are increasing, the $b_i$'s are
decreasing, the valencies $k_i$ are positive integers, and the multiplicities of the eigenvalues are positive integers. 
Besides the intersection arrays listed in the theorem, we only found the intersection arrays $\{ 5, 4,2; 1, 1, 2\}$
and $\{ 13, 12, 10; 1, 3, 4\}$. 
It was by Fon-der-Flaass \cite{F93} that there are no distance-regular graphs with intersection array $\{ 5, 4,2; 1, 
1, 2\}$. The intersection array  $\{ 13, 12, 10; 1, 3, 4\}$ is ruled out by \cite[Thm. 5.4.1]{bcn}.
This shows the theorem.
\epf.

\section{Proofs of Theorems 1.2 and 1.3}
In this section we give the proof of Theorems \ref{class} and \ref{3chrom}.\\
\\
{\bf Proof of Theorem \ref{class}:}
For diameter 2, it follows from Proposition \ref{diam=2}. For $a_1 \neq 0$ and diameter 3 and 4, it follows from 
Proposition \ref{a1}, and Theorems \ref{rnp} and \ref{a1class}.
 For diameter 3 and $a_1 =0$, it follows from Theorem \ref{classa=0}. 
 \epf

Before we give the proof of Theorem \ref{3chrom}, let us recall the chromatic number of a graph. 
A proper coloring with $t$ colors of a graph $\Gamma$ is a map $c: v(\Gamma) \rightarrow \{1, 2, \ldots, t\}$ where $t$ is a positive number  such that $c(x) \neq c(y)$ for any edge $xy$. The chromatic number  of $\Gamma$ denoted by $\chi(\Gamma)$ is the minimal $t$ such that there exists a proper coloring of $\Gamma$ with $t$ colors. We also say that such a graph is $\chi(\Gamma)$-chromatic. An independent set of $\Gamma$ is a set $S$ of vertices such that there are no edges between them.

Hoffman showed the following result for regular graphs.
\begin{lemma}{\em (Hoffman bound), cf. \cite[Prop. 1.3.2]{bcn}.}
Let $G$ be a $k$-regular graph with $n$ vertices and with smallest eigenvalue $\theta_{\min}$. Let $S$ be an 
independent set of
$\Gamma$ with $s$ vertices. Then $$s \leq \frac{n}{1 + \frac{k}{-\theta_{\min}}}.$$
\end{lemma}

This means that if a $k$-regular graph $\Gamma$ on $n$ vertices is 3-chromatic, then it must have an 
independent set of size $n/3$ and by the Hoffman bound we find that the smallest eigenvalue of $\Gamma$
is at most $-k/2$.  Now we are ready to give the proof for Theorem \ref{3chrom}. \\

{\bf Proof of Theorem \ref{3chrom}:} (i): By above we only need to 
check the graphs of Theorem \ref{class}.
The six graphs, we list, are shown to be 3-chromatic in \cite[Section 3.4]{BBH07}. For the case $a_1>0$, it was 
shown that the last three graphs are the only 3-chromatic distance-regular graphs with diameter 3 in \cite[Thm. 
3.6]{BBH07}. So we only need to check the graphs with $a_1 =0$. That the distance-regular graphs with 
intersection arrays $\{ 21, 20, 16; 1, 2, 12\}$ and $\{7, 6, 5; 1, 2, 3\}$ are not 3-chromatic follows from \cite[Sect. 
3.6]{BBH07}
That the distance-regular graph with intersection array $\{15, 14, 12; 1, 1, 9\}$ is not 3-chromatic follows from \cite[Sect. 3.7]{BBH07}. 

That the distance-regular graphs with intersection arrays  $\{5, 4, 2; 1, 1, 4\}$, $\{6, 5, 1; 1, 1, 6\}$, $\{7, 6, 6; 1, 1, 2\}$, $\{8, 7, 5; 1, 1, 4\}$ are not 3-chromatic follows from \cite[Sect. 3.9]{BBH07}.

(ii): The two graphs we list are shown be 3-chromatic in \cite[Sect. 3.4]{BBH07}. In \cite[Thm. 3.3 \& Prop. 3.8]{BBH07}, it is  shown that the Hamming graph $H(D, 3)$ is the only 3-chromatic distance-regular graph with  $c_2 \geq 2$, $a_D >0$ and $D \geq 4$. \cite[Thm. 3.3]{BBH07} also shows that the distance-regular graph with intersection array $\{6, 4, 2, 1; 1, 1, 4, 6\}$ is not 3-chromatic, as it has induced pentagons. That the regular near octagon associated with the Hall-Janko group (see \cite[Thm 13.6.1]{bcn}) with intersection array 
$\{ 10, 8, 8, 2; 1, 1, 4, 5\}$, is not 3-chromatic is shown on \cite[p. 299]{BBH07}. That a generalized octagon $GO(2, 4)$ with intersection array $\{10, 8, 8, 8; 1, 1, 1, 5\}$ is not 3-chromatic, follows from \cite[Thm. 3.2]{BBH07}.
This shows Theorem \ref{3chrom}.
\epf

\section{Open problems}
Now we give some open problems.
\begin{enumerate}
\item Classify the geometric distance-regular graphs with intersection number $a_1 =1$.
\item Finish the classification of the non-bipartite distance-regular graphs with diameter 4, valency $k$ and smallest eigenvalue at most $-k/2$.
\item Classify the non-bipartite distance-regular graphs with diameter 3, valency $k$ and smallest eigenvalue $-k/3$.
\end{enumerate}
\bibliographystyle{plain}
\bibliography{bib}

\end{document}